\documentclass{article}
\usepackage{amsthm,amsmath,amssymb}
\usepackage[left=0.5in,right=0.5in,bottom=1.0in,top=0.5in]{geometry} 
\usepackage{caption}
\usepackage{subcaption}
\usepackage{graphicx}
\usepackage{epstopdf}
\usepackage{color}
\usepackage{framed}
\usepackage{lscape}
\usepackage{rotating}
\usepackage{algorithm}
\usepackage[noend]{algpseudocode}
\usepackage{grffile}
\usepackage{multirow}
\usepackage{xcolor}
\usepackage{comment}
\usepackage{xr}
\usepackage[colorlinks=true,linkcolor=black,citecolor=black,urlcolor=black]{hyperref}

\newtheorem{thm}{Theorem}[section]

\newtheorem{lem}[thm]{Lemma}

\newtheorem{prop}[thm]{Proposition}

\theoremstyle{definition}
\newtheorem{rmk}{Remark}

\newcommand{\R}{\mathbb{R}}

\newcommand{\BLP}{\mathrm{BLP}}

\newcommand{\la}{\langle}
\newcommand{\argmin}{\operatorname*{argmin}}
\newcommand{\ra}{\rangle}

\newcommand{\op}{\mathrm{op}}

\renewcommand{\L}{\mathcal{L}}

\newcommand{\what}{\widehat}
\newcommand{\wtilde}{\widetilde}

\newcommand{\X}{\mathbf{X}}
\newcommand{\G}{\mathbf{G}}

\newcommand{\Y}{\mathbf{Y}}

\newcommand{\q}{\mathbf{q}}

\renewcommand{\u}{\mathbf{u}}
\renewcommand{\v}{\mathbf{v}}

\newcommand{\D}{\mathbf{D}}

\newcommand{\Fr}{\mathrm{F}}

\begin{document}
\title{Optimal singular value shrinkage for operator norm loss}
\author{William Leeb\thanks{School of Mathematics, University of Minnesota, Twin Cities. Minneapolis, MN.}}
\date{}
\maketitle

\abstract{
We study the denoising of low-rank matrices by singular value shrinkage. Recent work of Gavish and Donoho constructs a framework for finding optimal singular value shrinkers for a wide class of loss functions. We use this framework to derive the optimal shrinker for operator norm loss. The optimal shrinker matches the shrinker proposed by Gavish and Donoho in the special case of square matrices, but differs for all other aspect ratios. We precisely quantify the gain in accuracy from using the optimal shrinker. We also show that the optimal shrinker converges to the best linear predictor in the classical regime of aspect ratio zero.
}

\section{Introduction}

Low-rank matrix denoising is the task of estimating a low-rank matrix $\X$ from a noisy observed matrix $\Y = \X + \G$. In the setting of this paper, $\G$ is a matrix with iid Gaussian entries. We study the denoising procedure known as \emph{singular value shrinkage}, which keeps the singular vectors of $\Y$ while deflating the singular values to remove the effects of noise. Singular value shrinkage is a popular and well-studied methodology \cite{nadakuditi2014optshrink}, \cite{gavish2017optimal}, \cite{gavish2014optimal}, \cite{chatterjee2015matrix}, \cite{leeb2020matrix}, \cite{dobriban2020optimal}, \cite{leeb2019optimal}, \cite{shabalin2013reconstruction}, \cite{josse2016adaptive}, \cite{josse2016bootstrap}, \cite{bigot2017generalized}, \cite{donoho2014minimax}.

It has previously been observed, both for the matrix denoising problem and the related problem of low-rank covariance estimation, that the optimal singular value shrinker depends crucially on the choice of loss function between $\X$ and the estimated matrix $\what \X$ \cite{gavish2017optimal}, \cite{donoho2018optimal}. The work of Gavish and Donoho from \cite{gavish2017optimal} provides a general framework for deriving optimal singular value shrinkers for a wide class of loss functions. The results are applicable in a high-dimensional setting where the numbers of rows and columns grow to infinity, but the aspect ratio (the number of rows divided by the number of columns) converges to a definite limit. This is a standard setting, commonly referred to as the \emph{spiked model} \cite{johnstone2001distribution}.


In this paper, we revisit the special case of operator norm loss, defined as $\|\what \X - \X\|_{\op}$. We employ the the framework of Gavish and Donoho from \cite{gavish2017optimal} to derive the optimal singular value shrinker for this loss. We show that the optimal shrinker matches the shrinker proposed in \cite{gavish2017optimal} in the special case of square matrices, and provide a precise comparison of the two shrinkers for all aspect ratios. We will also show that when the columns of $\X$ are iid random vectors, then the optimal shrinker converges to the best linear predictor of each column in the limiting regime of aspect ratio zero, which can be interpreted as a ``classical'' statistical limit.

The rest of the manuscript is structured as follows. In Section \ref{sec-prelim}, we will formally state the spiked model assumptions, define singular value shrinkage, describe known results about the spiked model, and review the framework of Gavish and Donoho from \cite{gavish2017optimal}. In Section \ref{sec-shrinkage}, we will present the optimal shrinker for operator norm loss, and a comparison with the shrinker of \cite{gavish2017optimal}. In Section \ref{sec-blp}, we prove the convergence of optimal shrinkage to the best linear predictor in the classical regime. Section \ref{sec-proofs} contains detailed proofs of the main results. Section \ref{sec-conclusion} provides a brief conclusion.


\section{Preliminaries}
\label{sec-prelim}

\subsection{Model and estimation problem}

We observe a matrix $\Y = \X + \G$ of size $p$-by-$n$. The matrix $\G$ has entries which are iid $N(0,1/n)$. The matrix $\X$ is rank $r$, with singular value decomposition
\begin{align}
\X = \sum_{k=1}^{r} t_k \u_k \v_k^T.
\end{align}
Here, $\u_1,\dots,\u_r$ are orthonormal vectors in $\R^p$, $\v_1,\dots,\v_r$ are orthonormal vectors in $\R^n$, and $t_1 > \dots > t_r > 0$ are the singular values of $\X$.

An observation model of this form is generally referred to as a \emph{spiked model} \cite{johnstone2001distribution}. We study the spiked model in the asymptotic regime where both $n$ and $p = p(n)$ grow to infinity, and their ratio converges to a parameter $\gamma$, which we will refer to as the \emph{aspect ratio}:
\begin{align}
\gamma = \lim_{n \to \infty} \frac{p(n)}{n}.
\end{align}
Since $p$ and $n$ grow, for each $k=1,\dots,r$ we really have a sequence of singular vectors $\u_{k,n}$ and $\v_{k,n}$, indexed by $n$. However, to keep the notation to a minimum we will suppress the extra index $n$. Crucially, we assume that the singular values $t_1,\dots,t_r$ remain fixed, independently of $p$ and $n$.

Our goal is to estimate the low-rank matrix $\X$ from the noisy observed matrix $\Y$. We consider the use of operator norm loss, where the error between our estimator $\what \X$ and the true matrix $\X$ is given by:
\begin{align}
\L_{p,n}(\what \X, \X) = \|\what \X - \X\|_{\op}.
\end{align}
We will consider the class of singular value shrinkage estimators, which keep the top $r$ singular vectors of $\Y$ while changing the singular values. More precisely, we consider an estimator of the form
\begin{align}
\what \X^{\q} = \sum_{k=1}^{r} q_k \hat \u_k \hat \v_k^T,
\end{align}
where $\hat \u_1,\dots, \hat \u_r$ and $\hat \v_1,\dots, \hat \v_r$ are the top $r$ singular vectors of $\Y$, and $\q=(q_1,\dots,q_r)$ is the vector of singular values of $\what \X^{\q}$. As it turns out, for any specified choice of $\q = (q_1,\dots,q_r)$, the asymptotic loss
\begin{align}
\L_\infty(\q) = \lim_{p,n \to \infty} \L_{p,n}(\what \X^\q, \X)
\end{align}
is well-defined almost surely. The task is to find the values of $q_1,\dots,q_r$ that minimize the asymptotic loss; that is, we will find:
\begin{align}
\q^* = \argmin_{\q} \L_\infty(\q).
\end{align}
We also wish to evaluate the asymptotic loss $\L_\infty(\q^*)$ itself.

\subsection{Asymptotics of the spiked model}

The high-dimensional spiked model has been well-studied in the statistics and random matrix literature. It is known that there are precise relationships between the SVD of the observed matrix $\Y$ and the SVD of the low-rank matrix $\X$. If we write the SVD of $\Y$ as
\begin{align}
\Y = \sum_{k=1}^{\min(p,n)} \sigma_k \hat \u_k \hat \v_k^T,
\end{align}
then we can summarize the relevant results as follows:

\begin{prop}
\label{prop-vanilla}

For $1 \le k \le r$, the $k^{th}$ squared singular value of $\Y$ converges almost surely to the following deterministic limit:
\begin{align}
\label{eq-spikeforward}
\sigma_k^2= 
\begin{cases}
(t_k^2 + 1)\left( 1 + \frac{\gamma}{t_k^2}\right),
        & \text{ if } t_k  > \gamma^{1/4}, \\
(1 + \sqrt{\gamma})^2, & \text{ if }  t_k \le \gamma^{1/4},
\end{cases}
\end{align}

For $1 \le j, k \le r$, the squared cosines between the $j^{th}$ and $k^{th}$ singular vectors of $\X$ and $\Y$ converge almost surely to the following limits:

\begin{align}
\label{eq-c}
\lim_{p \to \infty}\la \hat \u_j, \u_k \ra^2=  c_{jk}^2 = 
\begin{cases}
\frac{1 - \gamma / t_k^4}{1 + \gamma / t_k^2},
        & \text{ if } j = k \text{ and } t_k  > \gamma^{1/4}, \\
0, & \text{ if } j \ne k \text{ or } t_k \le \gamma^{1/4},
\end{cases}
\end{align}
and
\begin{align}
\label{eq-ctilde}
\lim_{n \to \infty}\la \hat \v_j, \v_k \ra^2=\tilde c_{jk}^2 = 
\begin{cases}
\frac{1 - \gamma / t_k^4}{1 + 1 / t_k^2 },
        & \text{ if } j = k \text{ and } t_k > \gamma^{1/4}, \\
0, & \text{ if } j \ne k \text{ or } t_k \le \gamma^{1/4}.
\end{cases}
\end{align}

\end{prop}

A proof of this result may be found in \cite{paul2007asymptotics}, \cite{benaych2012singular}.

\begin{rmk}
\label{rmk-positive}
The signs of $c_{k}$ and $\tilde c_{k}$ are arbitrary, since the sign of a singular vector may be flipped. However, their product satisfies $c_{k} \tilde c_{k} \ge 0$, and we may therefore assume without loss of generality that $c_{k} \ge 0$ and $\tilde c_{k} \ge 0$ (see, e.g., \cite{nadakuditi2014optshrink}).
\end{rmk}

Proposition \ref{prop-vanilla} describes the behavior of the top $r$ singular components of $\Y$. For each singular value $\sigma_k$ of $\Y$ with $\sigma_k > 1 + \sqrt{\gamma}$, information about the corresponding component of $\X$ may be recovered. In particular, we may estimate $t_k$ by inverting formula \eqref{eq-spikeforward}:
\begin{align}
t_k = \sqrt{\frac{\sigma_k^2 - 1 - \gamma + \sqrt{(\sigma_k^2 - 1 - \gamma)^2 - 4\gamma}}{2}}.
\end{align}
From $t_k$, the cosines $c_k$ and $\tilde c_k$ are estimable by directly applying formulas \eqref{eq-c} and \eqref{eq-ctilde}.


\subsection{The shrinkage framework of \cite{gavish2017optimal}}

The work of Gavish and Donoho from \cite{gavish2017optimal} uses the behavior of the spiked model described in Proposition \ref{prop-vanilla} to construct a framework for finding asymptotically optimal singular value shrinkers. The key observation is that there are orthonormal bases of $\R^p$ and $\R^n$ in which the matrices $\what \X^{\q}$ and $\X$ may be expressed as follows:
\begin{align}
\X = \bigoplus_{k=1}^{r} 
\left(
\begin{array}{rr}
t_k  &   0  \\
0    &   0  \\
\end{array}
\right),
\end{align}
and
\begin{align}
\what \X^{\q} = \bigoplus_{k=1}^{r} 
q_k
\left(
\begin{array}{rr}
c_k \tilde c_k    &   c_k \tilde s_k  \\
s_k \tilde c_k    &   s_k \tilde s_k  \\
\end{array}
\right),
\end{align}
where $s = \sqrt{1 - c^2}$ and $\tilde s = \sqrt{1 - \tilde c^2}$. In other words, both $\X$ and $\what \X^{\q}$ are block-diagonal, with $r$ $2$-by-$2$ blocks of the prescribed form. The entries of these blocks depend only on the estimable quantities $t_k$, $c_k$, and $\tilde c_k$, and the values $q_k$ for which we solve. Because the operator norm is orthogonally-invariant and max-decomposable over block matrices, it follows that we may express the asymptotic operator norm loss as follows:
\begin{align}
\L_\infty(\q) = \max_{1 \le k \le r} 
\left\| 
\left(
\begin{array}{rr}
t_k  &   0  \\
0    &   0  \\
\end{array}
\right)
- 
q_k
\left(
\begin{array}{rr}
c_k \tilde c_k    &   c_k \tilde s_k  \\
s_k \tilde c_k    &   s_k \tilde s_k  \\
\end{array}
\right)
\right\|_{\op}.
\end{align}

Consequently, each optimal $q_k^*$ may in principle be solved for independently:
\begin{align}
\label{eq-argmin}
q_k^* = \argmin_{q_k} 
\left\| 
\left(
\begin{array}{rr}
t_k  &   0  \\
0    &   0  \\
\end{array}
\right)
- 
q_k
\left(
\begin{array}{rr}
c_k \tilde c_k    &   c_k \tilde s_k  \\
s_k \tilde c_k    &   s_k \tilde s_k  \\
\end{array}
\right)
\right\|_{\op}.
\end{align}
The remaining piece is to solve the minimization \eqref{eq-argmin}. In Section \ref{sec-shrinkage}, we will present a simple, estimable formula for the optimal $q_k^*$.

\section{Optimal shrinkage}
\label{sec-shrinkage}

\subsection{The optimal singular values}

In this section, we derive the optimal singular values $q_k^*$ and the resulting asymptotic operator norm loss $\L_\infty(\q^*)$. The main result is the following:

\begin{thm}
\label{thm-shrinker}
The optimal singular value shrinkage estimator $\what \X^{\q^*}$ of $\X$ from $\Y$ has singular values
\begin{align}
\label{eq-qstar}
q_k^* = 
\begin{cases}
t_k \sqrt{\frac{t_k^2 + \min\{1,\gamma\}}{t_k^2 + \max\{1,\gamma\}}}, 
    & \text{ if } \sigma_k > 1 + \sqrt{\gamma} \\
0, & \text{ otherwise }
\end{cases},
\end{align}
where $1 \le k \le r$. The loss $\|\what \X^{\q^*} - \X \|_{\op}$ converges almost surely to
\begin{align}
\L_\infty(\q^*) = \max_{1 \le k \le r} t_k \sqrt{1 - \min\{c_k^2,\tilde c_k^2\}}
=   t_1 \sqrt{1 - \min\{c_1^2,\tilde c_1^2\}} .
\end{align}

\end{thm}

We defer the proof of Theorem \ref{thm-shrinker} to Section \ref{proof-shrinker}.

\begin{figure}[h]
\centering
\includegraphics[scale=.5]{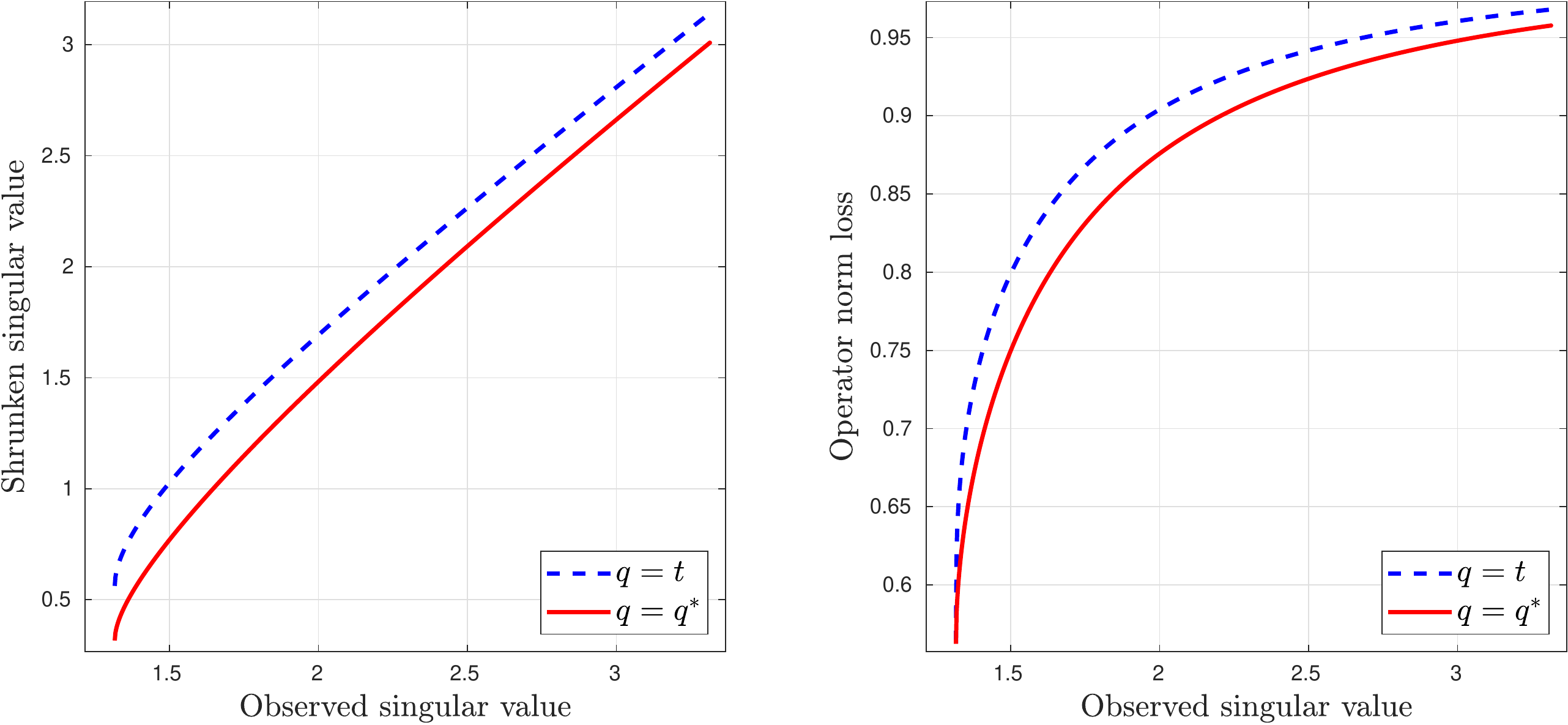}
\caption{Left: The singular values $q = q^*$ and $q = t$ as functions of the observed singular value $\sigma$. Right: The asymptotic operator norm losses as functions of $\sigma$.}
\label{fig-shrinkers}
\end{figure}

\begin{rmk}
The optimal $q_k^*$ and loss $\L_\infty(\q^*)$ may be consistently estimated from the observed singular values $\sigma_k$ of $\Y$ using \eqref{eq-spikeforward}, \eqref{eq-c}, and $\eqref{eq-ctilde}$.
\end{rmk}

\begin{rmk}
The paper \cite{gavish2017optimal} proposes the shrinker $\what \X^{\q}$ with singular values $q_k = t_k$. The optimal singular values $q_k^*$ we derive in Theorem \ref{thm-shrinker} are equal to $t_k$ when $\gamma = 1$ (the case of square matrices). In all other cases, however, the singular values $q_k = t_k$ are strictly suboptimal. In Figure \ref{fig-shrinkers}, we plot both shrinkers $q=q^*$ and $q=t$ and their asymptotic losses as functions as functions of the observed singular value $\sigma$.
\end{rmk}

\subsection{Comparison with the shrinker from \cite{gavish2017optimal}}

The paper \cite{gavish2017optimal} proposes use of the singular values $q_k = t_k$ for operator norm loss. This is a natural procedure, as it replace the noisy singular value $\sigma_k$ of $\Y$ with the ``true'' singular value $t_k$ of $\X$. However, according to Theorem \ref{thm-shrinker}, this choice is only optimal when $\gamma = 1$. We can quantify the gap in performance by comparing the relative errors. To simplify notation, we will consider only the rank $1$ setting, and subsequently drop subscripts; the same results hold when $\X$ is rank $r$.

\begin{prop}
\label{prop-compare}
Suppose $\X$ is a rank $1$ matrix with singular value $t > \gamma^{1/4}$. Let $\what \X^t$ denote the singular value shrinkage denoiser with $q=t$, and let $\what \X^{q^*}$ denote the optimal singular value shrinker from Theorem \ref{thm-shrinker}. Then
\begin{align}
\frac{\L_\infty(q^*)}{\L_\infty(t)}
= \lim_{n \to \infty} \frac{\|\X - \what \X^{q^*}\|_{\op}}{\|\X - \what \X^{t}\|_{\op}}
= \sqrt{\frac{1 + \min\{c,\tilde c\}}{ 1 + \max\{c,\tilde c\}} },
\end{align}
where the limit holds almost surely as $p/n \to \gamma$.
\end{prop}

In the next result, we derive a limiting expression for the error ratio as $\gamma \to 0$.

\begin{prop}
\label{prop-compare-classical}
Define the asymptotic error ratio $R$:
\begin{align}
R = \sqrt{\frac{1 + \min\{c,\tilde c\}}{ 1 + \max\{c,\tilde c\}} }.
\end{align}
For fixed $t > 0$, $R = R(\gamma,t)$ is an increasing function of $\gamma < \min\{t^4,1\}$, and its minimum value is
\begin{align}
\lim_{\gamma \to 0} R(\gamma,t) = \sqrt{\frac{1}{2} \left(1 + \sqrt{\frac{t^2}{t^2 + 1}} \right)}.
\end{align}
\end{prop}

\begin{rmk}
If the columns of $\X$ are $n$ iid random vectors from a rank $r$ distribution in $\R^p$, then the $\gamma = 0$ limit considered in Proposition \ref{prop-compare-classical} can be thought of informally as the ``classical'' setting, where the number of observations $n$ grows faster than the number of features $p$.
\end{rmk}

\begin{figure}[h]
\centering
\includegraphics[scale=.5]{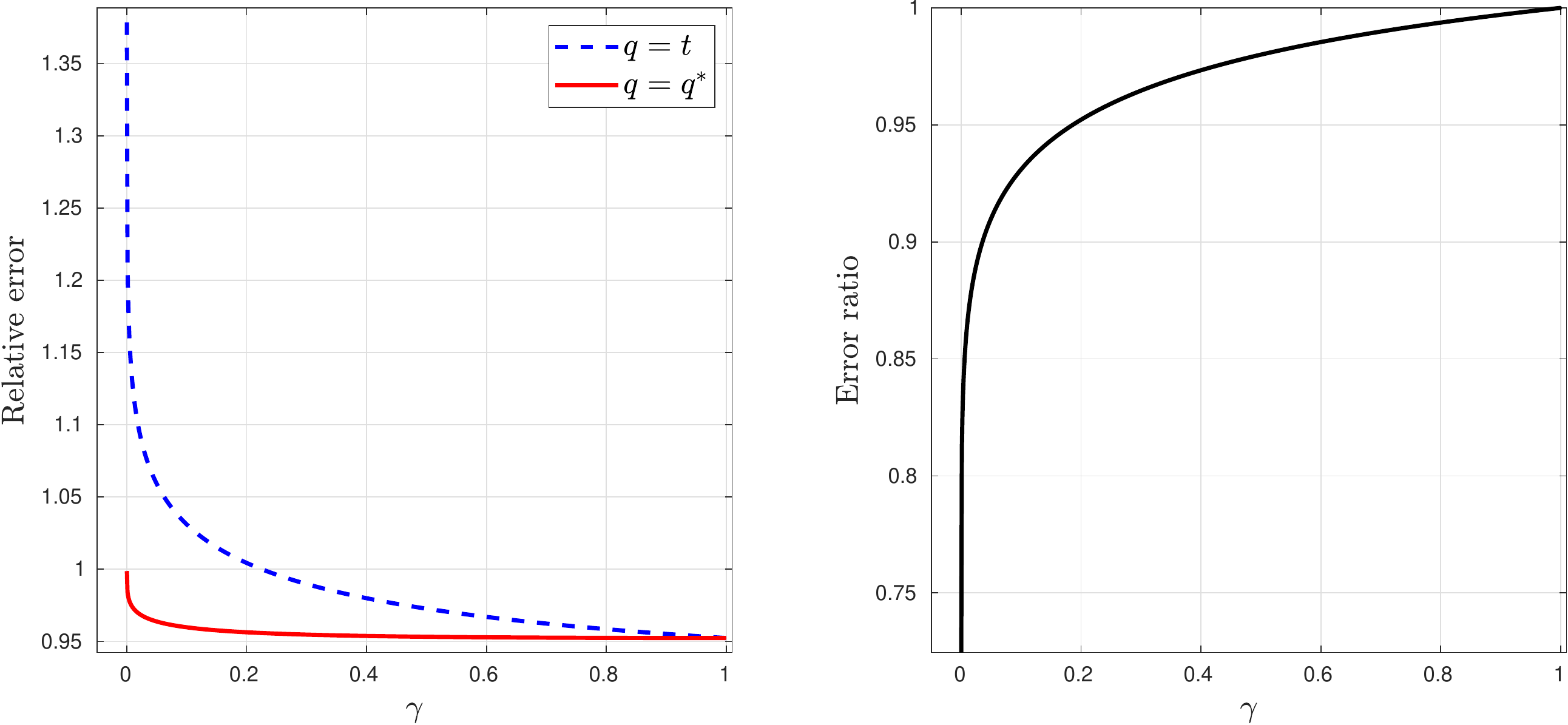}
\caption{Left: Relative errors $\|\what \X^q - \X\|_{\op} / t$ for $q = q^*$ and $q=t$ as a function of $\gamma$, where $t = \gamma^{1/4} + 1/20$. Right: The ratio of the errors $\|\what \X^{q^*} - \X\|_{\op} / \|\what \X^{t} - \X\|_{\op}$ as a function of $\gamma$.}
\label{fig-ratio}
\end{figure}

Proposition \ref{prop-compare-classical} shows that the performance of the optimal shrinker $q=q^*$ over $q=t$ should be most evident when $\gamma$ and $t$ are both small. To illustrate this, in Figure \ref{fig-ratio} we compare error curves of $q = t$ and $q = q^*$ as a function of $\gamma \in [0,1]$, where the signal strength $t$ is chosen to decrease with $\gamma$ as $t = \gamma^{1/4} + 1/20$. The left panel shows the relative errors (the error divided by the signal norm $t$) for the two shrinkers as a function of $\gamma$, and the right panel shows the ratio of the errors. As $\gamma$ and $t$ both decrease, the relative performance of the optimal shrinker $q=q^*$ increases over the shrinker $q=t$.

\section{Convergence to the best linear predictor}
\label{sec-blp}

In this section, we consider the setting where the columns of $\X$ are iid random vectors from a distribution in $\R^p$ with mean zero. We will write each column of $\X$ as $X_j / \sqrt{n}$, where $X_j$ is a random vector of the following form:
\begin{align}
X_j = \sum_{k=1}^{r} t_k z_{jk} \u_k,
\end{align}
where the $z_{jk}$ are mean zero, unit variance sub-Gaussian random variables, and the $\u_k$ are the orthonormal principal components of $X_j$.

\begin{rmk}
In the new setting, each $t_k$ is the standard deviation of $X_j$ along the principal component $\u_k$, not the singular value of $\X$; and $\u_1,\dots,\u_r$ are not the left singular vectors of $\X$. However, in the large $n$ limit, the singular values of $\X$ converge almost surely to $t_1,\dots,t_r$, and the left singular vectors of $\X$ almost surely make zero angle with, respectively, $\u_1,\dots,\u_r$. In this sense, the new notation is consistent with the old notation.
\end{rmk}

\begin{rmk}
The assumption that $X_j$ has mean zero is easily removed by subtracting the sample mean from each $X_j$.
\end{rmk}

It is known \cite{mackay2004deconv} that the best linear predictor of $X_j$ from $Y_j$ has the following form:
\begin{align}
\label{eq-blp}
\what X_j^{\BLP} = \sum_{k=1}^{r} \frac{t_k^2}{t_k^2 + 1} \la Y_j , \u_k\ra \u_k. 
\end{align}
The next result shows that optimal singular value shrinkage with operator norm loss converges to the best linear predictor when $\gamma = 0$. We will let $\what \X^{\q^*} = [\what X_j^{\q^*},\dots, \what X_j^{\q^*}] / \sqrt{n}$ and $\what \X^{\BLP} = [\what X_j^{\BLP},\dots, \what X_j^{\BLP}] / \sqrt{n}$.

\begin{thm}
\label{thm-blp}
In the limit $n \to \infty$ and $p/n \to 0$,
\begin{align}
\lim_{n \to \infty} \frac{1}{n} \sum_{j=1}^{n} \|\what X_j^{\q^*} - \what X_j^{\BLP}\|^2
= \lim_{n \to \infty} \|\what \X_j^{\q^*} - \what \X_j^{\BLP}\|_{\Fr}^2
= 0,
\end{align}
where the limit holds almost surely.
\end{thm}

Theorem \ref{thm-blp} is proven in Section \ref{proof-blp}. It is a consequence of the following result, whose proof follows from \cite{leeb2019optimal}:
\begin{lem}
\label{lem-linear}
Let $\what \X^{\q}$ be any singular value shrinker, with singular values $q_1,\dots,q_r$. Define the linear predictor $\wtilde X_j^{\q}$ by:
\begin{align}
\label{eq-linear}
\wtilde X_j^{\q} = \sum_{k=1}^{r} \frac{q_k}{\sigma_k} \la Y_j , \u_k\ra \u_k,
\end{align}
where $\sigma_1,\dots,\sigma_r$ are the top $r$ singular values of $\Y$, and $\hat \u_1,\dots, \hat \u_r$ are the top $r$ left singular vectors of $\Y$. Then
\begin{align}
\lim_{n \to \infty} \frac{1}{n} \sum_{j=1}^{n} \|\what X_j^{\q} - \wtilde X_j^{\q}\|^2
= \lim_{n \to \infty} \|\what \X_j^{\q} - \wtilde \X_j^{\q}\|_{\Fr}^2
= 0,
\end{align}
where the limit holds almost surely in the limit $n \to \infty$, $p/n \to 0$.
\end{lem}

\begin{rmk}
\label{rmk-blp}
Lemma \ref{lem-linear} states that singular value shrinkage converges to a linear predictor of the form \eqref{eq-linear} when $\gamma = 0$ (the ``classical'' regime of zero aspect ratio). Comparing \eqref{eq-linear} to the form of the best linear predictor \eqref{eq-blp}, we see that the shrinker $\what \X^{\q}$ will converge to the best linear predictor if and only if
\begin{align}
q_k = \frac{ \sigma_k t_k^2}{t_k^2 + 1} = t_k \sqrt{\frac{t_k^2}{t_k^2 + 1}}, \quad 1 \le k \le r,
\end{align}
when $\gamma = 0$. Any other choice of shrinker, including $q_k = t_k$, will result in convergence to a suboptimal linear filter in the $\gamma = 0$ regime.
\end{rmk}

\begin{figure}[h]
\centering
\includegraphics[scale=.5]{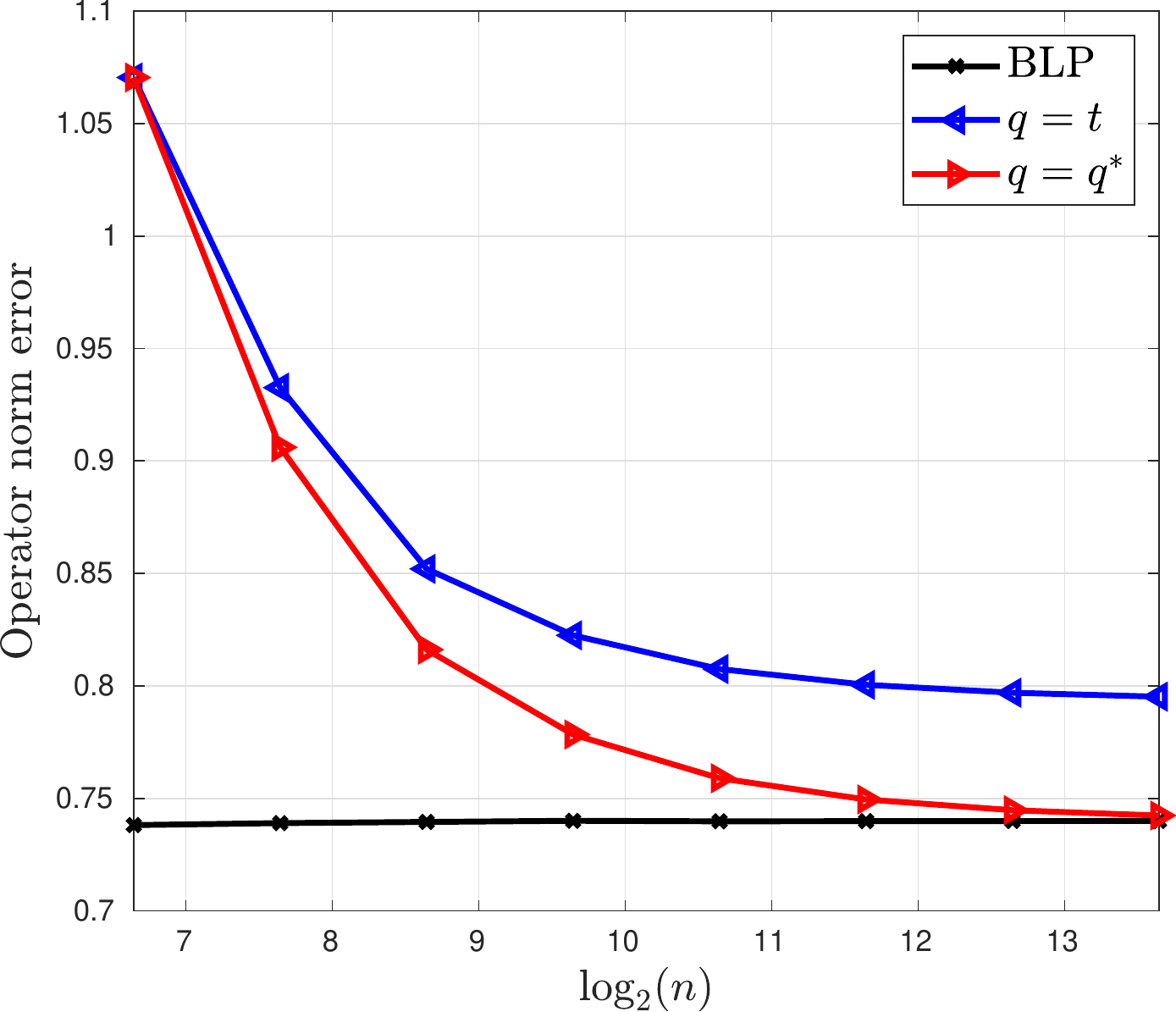}
\caption{Operator norm errors for $p = 100$ and increasing $n$ for the BLP \eqref{eq-blp} and the shrinkers $\what \X^{q}$, $q=q^*$ and $q=t$. The signal is rank $1$ with $t=1.1$. Errors are averaged over $4000$ runs.}
\label{fig-blp}
\end{figure}

To illustrate Theorem \ref{thm-blp} numerically, we draw $n$ iid observations from a spiked model in $\R^p$, for increasing values of $n \ge p$. We take the rank $r=1$ and $t=t_1 = 1.1$ (to ensure the signal is detectable for all $n \ge p$). We apply the best linear predictor $\what X_j^{\BLP}$ (which assumes the principal component $\u = \u_1$ is known), optimal shrinkage $\what X_j^{q^*}$, and the suboptimal shrinker $\what X_j^{t}$. In Figure \ref{fig-blp}, we plot the average operator norm error over 4000 runs of the experiment, as a function of $n$.

The error for the BLP is approximately constant, since the BLP does not vary with the sample size $n$. As $n$ grows, the error for optimal shrinkage approaches that of the BLP, because $\what X_j^{q^*}$ converges to the BLP $\what X_j^{\BLP}$. By contrast, the error for the suboptimal shrinker converges to a strictly larger value, since $\what X_j^{t}$ converges to the suboptimal linear predictor $\wtilde X_j^{t}$.

\section{Proofs}
\label{sec-proofs}

\subsection{Proof of Theorem \ref{thm-shrinker}}
\label{proof-shrinker}

Theorem \ref{thm-shrinker} is a consequence of the following result:

\begin{lem}
\label{lem-opmin}
Let $t > 0$, and $0 \le c \le 1$, $0 \le \tilde c \le 1$. For any $q \in \R$, define the $2$-by-$2$ matrix $\D(q)$ by
\begin{align}
\D(q) = 
\left( 
\begin{array}{rr}
t  & 0  \\
0  & 0
\end{array}
\right)
-
q
\left( 
\begin{array}{rr}
c \tilde c  & c \tilde s  \\
s \tilde c  & s \tilde s
\end{array}
\right),
\end{align}
where $s = \sqrt{1 - c^2}$ and $\tilde s = \sqrt{1 - \tilde c^2}$. Define the function $F$ by:
\begin{align}
\label{eq-F}
F(q) = \|\D(q)\|_{\op}^2,
\end{align}
the squared operator norm of $\D(q)$. Then the value $q^*$ that minimizes $F$ defined in \eqref{eq-F} is given by 
\begin{align}
q^* = t \cdot \frac{\min\{c,\tilde c\}}{\max\{c,\tilde c\}},
\end{align}
when $c$ and $\tilde c$ are not both $0$; and $q^*$ may be taken as any value with $|q^*| \le t$ if $c =  \tilde c = 0$. Furthermore, the value of $F$ at $q^*$ is given by:
\begin{align}
\label{eq-error}
F(q^*) = t^2 \cdot \max\{s^2,\tilde s^2\}.
\end{align}
\end{lem}

\begin{rmk}
Lemma \ref{lem-opmin} is applicable for general parameters $t > 0$ and $0 \le c \le 1$, $0 \le \tilde c \le 1$, even when they do not satisfy the relationships \eqref{eq-spikeforward}, \eqref{eq-c} and \eqref{eq-ctilde}.
\end{rmk}

We will first prove Theorem \ref{thm-shrinker} assuming Lemma \ref{lem-opmin}, and then return to the proof of Lemma \ref{lem-opmin} itself in Section \ref{proof-opmin}. Without loss of generality, we assume $\gamma \le 1$. If $\sigma_k > 1 + \sqrt{\gamma}$, or equivalently $t_k > \gamma^{1/4}$, Lemma \ref{lem-opmin} says that the optimal singular value is given by $q_k^* = t_k \tilde c_k / c_k$. Using the expressions \eqref{eq-c} and \eqref{eq-ctilde} for $c_k$ and $\tilde c_k$, formula \eqref{eq-qstar} follows immediately. When $\sigma \le 1+\sqrt{\gamma}$, or equivalently $t_k \le \gamma^{1/4}$, both $c_k$ and $\tilde c_k$ are zero. Consequently, $q_k^* = 0$ is optimal.

All the remains is to show the error formula:
\begin{align}
\max_{1 \le k \le r} t_k \sqrt{1 - \min\{c_k^2,\tilde c_k^2\}}
=   t_1 \sqrt{1 - \min\{c_1^2,\tilde c_1^2\}} .
\end{align}
Without loss of generality, assume $\gamma \le 1$, so $\tilde c \le c$. Then
\begin{align}
t^2 (1 - \tilde c^2)
= t^2 \frac{t^2 + \gamma}{t^4 + t^2}
= \frac{t^2 + \gamma}{t^2 + 1}
\end{align}
is an increasing function of $t$, and is equal to $\sqrt{\gamma}$ (the largest possible error for any component $t_k \le \gamma^{1/4}$) when $t = \gamma^{1/4}$. Consequently, the maximum error is achieved at $t = t_1$.

\subsubsection{Proof of Lemma \ref{lem-opmin}}
\label{proof-opmin}

First, suppose $c = \tilde c = 0$. Then
\begin{align}
\D(q) = 
\left( 
\begin{array}{rr}
t  & 0  \\
0  & -q
\end{array}
\right),
\end{align}
and so $F(q) = \|\D(\q)\|_{\op}^2 = \max\{t^2,q^2\}$. Consequently, any $q$ with $|q|  \le t$ minimizes $F(q)$, and since $s = \tilde s = 1$, $F(q) = t^2 \cdot \max\{s^2 \tilde s^2\}$ for such $q$.

Next, assume that $c$ and $\tilde c$ are not both $0$. The proof when $c = \tilde c$ is identical to the proof for optimal shrinkage of eigenvalues for covariance estimation contained in \cite{donoho2018optimal}; so we will assume that $c \ne \tilde c$. We will use the expression for the operator norm as a function of $q$ derived in \cite{gavish2017optimal}. We have:
\begin{align}
F(q) &= 
\frac{q^2 + t^2 - 2qtc\tilde c + \sqrt{(q^2 + t^2 - 2qtc\tilde c)^2 - 4(t q s \tilde s)^2}}{2}
\nonumber \\
& =\frac{A(q) + \sqrt{A(q)^2 - 4B(q)^2}}{2},
\end{align}
where 
\begin{align}
A(q) = q^2 + t^2 - 2qtc\tilde c,
\end{align}
and 
\begin{align}
B(q) = -t q s \tilde s.
\end{align}
The function $F(q)$ is differentiable whenever $A(q)^2 - 4B(q)^2 > 0$. Furthermore, at points where $F$ is differentiable, its derivative is given by
\begin{align}
F'(q) = \frac{1}{2} \left(A'(q)
    + \frac{A(q) A'(q) - 4 B(q) B'(q)}{\sqrt{A(q)^2 - 4B(q)^2}} \right).
\end{align}
%

Suppose, without loss of generality, that $\tilde c < c$. Then $q^* = t \tilde c / c$. First, we have
\begin{align}
A'(q) = 2q - 2tc\tilde c.
\end{align}
Consequently
\begin{align}
A'(q^*) = 2 t \frac{\tilde c}{ c} - 2 t c \tilde c
= 2 t \frac{\tilde c}{c}\left( 1 - c^2  \right)
= 2 t \frac{\tilde c}{c} s^2.
\end{align}
We also have
\begin{align}
\label{eq1000}
A(q^*) = t^2 \left( 1 + \frac{\tilde c^2}{c^2} - 2 \tilde c^2\right),
\end{align}
and so
\begin{align}
\label{eq1100}
A(q^*) A'(q^*) 
= 2 t^3 s^2 \frac{\tilde c}{c} \left( 1 + \frac{\tilde c^2}{c^2} - 2 \tilde c^2 \right).
\end{align}

Next, observe that for all $q$,
\begin{align}
B'(q) = - ts \tilde s.
\end{align}
We also have
\begin{align}
B(q^*) = -t^2 \frac{\tilde c}{c} s \tilde s,
\end{align}
and so
\begin{align}
\label{eq1200}
B(q^*) B'(q^*) = t^3 s^2 \frac{\tilde c}{c} (1-\tilde c^2).
\end{align}

Combining \eqref{eq1100} and \eqref{eq1200}, we obtain:
\begin{align}
\label{eq1400}
A(q^*)A'(q^*) - 4B(q^*)B'(q^*)
&= 2 t^3 s^2 \frac{\tilde c}{c} \left( 
1 + \frac{\tilde c^2}{c^2} - 2 \tilde c^2 - 2 + 2 \tilde c^2 \right)
\nonumber \\
&= 2 t^3 s^2 \frac{\tilde c}{c} \left( 
\frac{\tilde c^2}{c^2} - 1 \right)
\nonumber \\
&= 2 t^3 s^2 \frac{\tilde c}{c^3} (\tilde c^2 - c^2).
\end{align}

Next, we observe that we may write
\begin{align}
A(q^*)^2 - 4B(q^*)^2 = t^4\left(1 - \frac{\tilde c^2}{c^2}\right)^2,
\end{align}
and consequently,
\begin{align}
\label{eq1600}
\sqrt{A(q^*)^2 - 4B(q^*)^2}
= \frac{t^2}{c^2}(c^2 - \tilde c^2),
\end{align}
where we have used the fact that $c > \tilde c$. Note that \eqref{eq1600} implies that $F$ is differentiable at $q^*$ whenever $c \ne \tilde c$.

Combining \eqref{eq1400} and \eqref{eq1600}, we get:
\begin{align}
\frac{A(q^*)A'(q^*) - 4B(q^*)B'(q^*)}{\sqrt{A(q^*)^2 - 4B(q^*)^2}}
= -2t\frac{\tilde c}{c}s^2.
\end{align}

Consequently,
\begin{align}
F'(q^*) &= \frac{1}{2} \left(A'(q^*)
    + \frac{A(q^*) A'(q^*) - 4 B(q^*) B'(q^*)}{\sqrt{A(q^*)^2 - 4B(q^*)^2}} \right)
\nonumber \\
&= \frac{1}{2} \left( 2 t \frac{\tilde c}{c} s^2
     -2t\frac{\tilde c}{c} s^2 \right)
\nonumber \\
&=0.
\end{align}

Next, we evaluate $F(q^*)$. From \eqref{eq1000} and \eqref{eq1600}, we get
\begin{align}
F(q^*) &= \frac{1}{2} \left(t^2 \left( 1 + \frac{\tilde c^2}{c^2} - 2 \tilde c^2\right)  
    + \frac{t^2}{c^2}(c^2 - \tilde c^2)\right)
\nonumber \\
&= \frac{1}{2} \left(t^2 + t^2\frac{\tilde c^2}{c^2} - 2 t^2\tilde c^2
    + t^2  - t^2 \frac{\tilde c^2}{c^2}\right)
\nonumber \\
&= t^2 (1 - \tilde c^2),
\end{align}
which is the desired result.

\subsection{Proof of Propositions \ref{prop-compare} and \ref{prop-compare-classical}}
\subsubsection{Proof of Proposition \ref{prop-compare}}
The error formula for $q = t$ found in \cite{gavish2017optimal}:
\begin{align}
\L_\infty(t) = t \cdot \sqrt{1 - c \tilde c + |c - \tilde c|}.
\end{align}
Assuming, without loss of generality, that $\gamma \le 1$ and hence $\tilde c \le c$, then
\begin{align}
\L_\infty(t) = t \cdot \sqrt{1 - c \tilde c + c - \tilde c}
= t \cdot \sqrt{(1-\tilde c)(1+c)}.
\end{align}
Since $\L_\infty(q^*) = t \sqrt{1 - \tilde c^2} = t \sqrt{(1 - \tilde c)(1 + \tilde c)}$ from Theorem \ref{thm-shrinker}, dividing $\L_\infty(t)$ by $\L_\infty(q^*)$ gives the result.

\subsubsection{Proof of Proposition \ref{prop-compare-classical}}

The value $t > 0$ is fixed, and we consider the range $\gamma < \min\{t^4,1\}$. Define the function
\begin{align}
E(\gamma) = R(\gamma,t)^2 = \frac{1 + \tilde c (\gamma)}{ 1 + c(\gamma)},
\end{align}
where we treat the cosines $c = \sqrt{(t^4 - \gamma) / (t^4 + \gamma t^2)}$ and $\tilde c = \sqrt{(t^4 - \gamma) / (t^4 + t^2)}$ as functions of $\gamma$. Then the derivative of $E(\gamma)$ is:
\begin{align}
E'(\gamma) 
= \frac{1}{2 c(\gamma) \tilde c(\gamma) (1 + c(\gamma))^2}
\left(  \frac{\tilde c(\gamma)(1+\tilde c(\gamma))(t^2+1)}{(t^2 + \gamma)^2 }
    - \frac{c(\gamma)(1+c(\gamma))}{t^4 + t^2} \right).
\end{align}
Consequently, to show $E'(\gamma) > 0$, it is enough to show that
\begin{align}
\frac{\tilde c(\gamma)(1+\tilde c(\gamma))(t^2+1)}{(t^2 + \gamma)^2 }
    - \frac{c(\gamma)(1+c(\gamma))}{t^4 + t^2} > 0,
\end{align}
or equivalently that
\begin{align}
\frac{\tilde c(\gamma)}{ c(\gamma)} 
\frac{1 + \tilde c(\gamma)}{1 + c(\gamma)} 
\left(\frac{t^2 + 1}{t^2 + \gamma}\right)^2
t^2 - 1 > 0.
\end{align}
Since $\gamma < 1$, $\tilde c(\gamma) < c(\gamma)$; consequently,
\begin{align}
\frac{\tilde c(\gamma)}{ c(\gamma)} \frac{1 + \tilde c(\gamma)}{1 + c(\gamma)}  
>  \frac{\tilde c(\gamma)^2}{ c(\gamma)^2} 
= \frac{t^2 + \gamma}{t^2 + 1}.
\end{align}
Hence it is enough to show
\begin{align}
\left(\frac{t^2 + 1}{t^2 + \gamma}\right)
t^2 - 1 > 0;
\end{align}
but this follows immediately, since $\gamma < t^4$.

The limit of $R(\gamma,t)$ as $\gamma \to 0$ follows immediately from the limits:
\begin{align}
\lim_{\gamma \to 0} c(\gamma) = \lim_{\gamma \to 0} \sqrt{\frac{t^4 - \gamma}{t^4 + \gamma t^2}}
= 1
\end{align}
and
\begin{align}
\lim_{\gamma \to 0} \tilde c(\gamma) = \lim_{\gamma \to 0} \sqrt{\frac{t^4 - \gamma}{t^4 + t^2}}
= \sqrt{\frac{t^2}{t^2 + 1}}.
\end{align}

\subsection{Proof of Theorem \ref{thm-blp}}
\label{proof-blp}

Since each component is treated separately, we drop the subscripts. Following Remark \ref{rmk-blp}, Theorem \ref{thm-blp} will be proven if we show that:
\begin{align}
q^* = \frac{ \sigma t^2}{t^2 + 1} = t \sqrt{\frac{t^2}{t^2 + 1}}
\end{align}
when $\gamma = 0$. To prove this identity, observe that when $\gamma = 0$, $c = 1$, and
\begin{align}
\tilde c = \sqrt{\frac{t^2}{t^2 + 1}}.
\end{align}
Since $q^* = t \cdot \min\{c, \tilde c\} = t \tilde c$, the result follows immediately.

\section{Conclusion}
\label{sec-conclusion}

We have considered the problem of estimating a low-rank matrix $\X$ from noisy observations $\Y = \X + \G$, where we measure the error by operator norm loss $\|\X - \what \X\|_{\op}$. We have proven (Theorem \ref{thm-shrinker}) that the optimal singular value shrinker has singular values of the form \eqref{eq-qstar}. For square matrices ($\gamma = 1$), the optimal singular values agree with those proposed in \cite{gavish2017optimal}, though the two methods and the resulting errors differ increasingly as $\gamma$ differs from $1$, or equivalently as the matrix becomes more rectangular (Proposition \ref{prop-compare-classical}). We have also shown  (Theorem \ref{thm-blp}) that when the columns of $\X$ are iid vectors in $\R^p$, then in the classical regime ($\gamma=0$) the optimal shrinker converges to the best linear predictor.

\section*{Acknowledgements}

I acknowledge support from the NSF BIGDATA program IIS 1837992 and BSF award 2018230.

\bibliographystyle{plain}
\bibliography{refs}

\begin{thebibliography}{10}

\bibitem{benaych2012singular}
Florent Benaych-Georges and Raj~Rao Nadakuditi.
\newblock The singular values and vectors of low rank perturbations of large
  rectangular random matrices.
\newblock {\em Journal of Multivariate Analysis}, 111:120--135, 2012.

\bibitem{bigot2017generalized}
J\'er\'emie Bigot, Charles Deledalle, and Delphine F\'eral.
\newblock Generalized {SURE} for optimal shrinkage of singular values in
  low-rank matrix denoising.
\newblock {\em Journal of Machine Learning Research}, 18:1--50, 2017.

\bibitem{chatterjee2015matrix}
Sourav Chatterjee.
\newblock Matrix estimation by universal singular value thresholding.
\newblock {\em The Annals of Statistics}, 43(1):177--214, 2015.

\bibitem{dobriban2020optimal}
Edgar Dobriban, William Leeb, and Amit Singer.
\newblock Optimal prediction in the linearly transformed spiked model.
\newblock {\em Annals of Statistics}, 48(1):491--513, 2020.

\bibitem{donoho2018optimal}
David~L. Donoho, Matan Gavish, and Iain~M. Johnstone.
\newblock Optimal shrinkage of eigenvalues in the spiked covariance model.
\newblock {\em Annals of Statistics}, 46(4):1742–1778, 2018.

\bibitem{donoho2014minimax}
Matan Gavish and David~L. Donoho.
\newblock Minimax risk of matrix denoising by singular value thresholding.
\newblock {\em The Annals of Statistics}, 42(6):2413--2440, 2014.

\bibitem{gavish2014optimal}
Matan Gavish and David~L. Donoho.
\newblock The optimal hard threshold for singular values is $4/\sqrt {3} $.
\newblock {\em IEEE Transactions on Information Theory}, 60(8):5040--5053,
  2014.

\bibitem{gavish2017optimal}
Matan Gavish and David~L. Donoho.
\newblock Optimal shrinkage of singular values.
\newblock {\em IEEE Transactions on Information Theory}, 63(4):2137--2152,
  2017.

\bibitem{johnstone2001distribution}
Iain~M Johnstone.
\newblock On the distribution of the largest eigenvalue in principal components
  analysis.
\newblock {\em Annals of Statistics}, 29(2):295--327, 2001.

\bibitem{josse2016adaptive}
Julie Josse and Sylvain Sardy.
\newblock Adaptive shrinkage of singular values.
\newblock {\em Statistics and Computing}, 26:715–724, 2016.

\bibitem{josse2016bootstrap}
Julie Josse and Stefan Wager.
\newblock Bootstrap-based regularization for low-rank matrix estimation.
\newblock {\em Journal of Machine Learning Research}, 17:1--29, 2016.

\bibitem{leeb2020matrix}
William Leeb.
\newblock Matrix denoising for weighted loss functions and heterogeneous
  signals.
\newblock {\em arXiv:1902.09474}, 2020.

\bibitem{leeb2019optimal}
William Leeb and Elad Romanov.
\newblock Optimal spectral shrinkage and {PCA} with heteroscedastic noise.
\newblock {\em arXiv:1811.02201}, 2019.

\bibitem{mackay2004deconv}
D.~J.~C. MacKay.
\newblock Deconvolution.
\newblock In {\em Information Theory, Inference and Learning Algorithms}, pages
  550--551. Cambridge University Press, Camridge, UK, 2004.

\bibitem{nadakuditi2014optshrink}
Raj~Rao Nadakuditi.
\newblock {OptShrink}: An algorithm for improved low-rank signal matrix
  denoising by optimal, data-driven singular value shrinkage.
\newblock {\em IEEE Transactions on Information Theory}, 60(5):3002--3018,
  2014.

\bibitem{paul2007asymptotics}
Debashis Paul.
\newblock Asymptotics of sample eigenstructure for a large dimensional spiked
  covariance model.
\newblock {\em Statistica Sinica}, 17(4):1617--1642, 2007.

\bibitem{shabalin2013reconstruction}
Andrey~A. Shabalin and Andrew~B. Nobel.
\newblock Reconstruction of a low-rank matrix in the presence of {G}aussian
  noise.
\newblock {\em Journal of Multivariate Analysis}, 118:67--76, 2013.

\end{thebibliography}

\end{document}